\documentclass[preprint,12pt]{elsarticle}

\usepackage{graphics}
\usepackage{graphicx}
\usepackage{epsfig}

\usepackage{amssymb}
\usepackage{amsmath}
\newtheorem{thm}{Theorem}
\journal{Journal of Computational and Applied Mathematics}

\begin{document}

\begin{frontmatter}



\title{Discontinuous Galerkin Method for Total Variation Minimization on one-dimensional Inpainting Problem}


\author{Xijian Wang}

\address{School of Mathematics and Computing Science, Wuyi University, People's Republic of China}
\ead{wangxj1980426@gmail.com}
\begin{abstract}
This paper is concerned with the numerical minimization of energy functionals in $BV(\Omega)$ (the space of bounded variation functions) involving total variation for gray-scale 1-dimensional inpainting problem. Applications are shown by finite element method and discontinuous Galerkin method for total variation minimization. We include the numerical examples which show the different recovery image by these two methods.
\end{abstract}


\end{frontmatter}

% \linenumbers




%\maketitle
%% -----------------------------------------------------------
%
\section{Introduction}
In the first chapter of the book ~\cite{For10} Holger Rauhut has already introduced that the minimization of $\ell_{1}$-norms occupies a fundamental role for the promotion of sparse solutions. This understanding furnishes an important interpretation of total variation minimization ~\cite{Rud92} as a regularization technique for image inpainting.
In this paper we consider as in ~\cite{Cha97, Ves01} the minimization in $BV(\Omega)$ (the space of bounded variation functions ~\cite{Aub06,Chan05}) of the functional
\begin{equation}\label{original}
\mathcal{J}(u):=\int_{\Omega}\left\vert Tu(x)-g(x)\right\vert^{2}dx+2\lambda\left\vert Du\right\vert(\Omega),
\end{equation}
where $\Omega\subset R^{d}$, for $d=1,2$ be a bounded Lipschitz domain, $T:L^{2}(\Omega)\rightarrow L^{2}(\Omega)$ is a bounded linear operator , $g\in L^{2}(\Omega)$ is a datum, $\left\vert Du\right\vert(\Omega):=\int_{\Omega}\left\vert \nabla u(x)\right\vert dx$ is the total variation of $u$, and $\lambda >0$ is a fixed regularization parameter~\cite{Eng96}. Several numerical strategies to efficiently perform total variation minimization have been proposed in the literature, refer to ~\cite{Cham97,Gol09,Osh05,Wei09}. However, the interesting solutions may be discontinuous, e.g., along curves in 2D. Hence, the crucial difficulty is the correct numerical treatment of interfaces, with the preservation of crossing discontinuities and the correct matching where the solution is continuous instead, see Section 7.1.1 in \cite{For09}. In order to deal promptly with the discontinuity, we have studied the applications to gray-scale 1-dimensional inpainting problem by the finite element method and discontinuous Galerkin method(Refer to \cite{Coc00,Riv08}) for total variation minimization, respectively.

The paper is organized as follows. Section 2 is devoted to the alternating-minimization algorithm to compute minimizers of $\mathcal{J}(u)$. In Section 3, the method of finite element method for total variation minimization for our problem is illustrated. The main work of this paper about discontinuous Galerkin method for total variation minimization is described in Section 4. Finally, we include some numerical experiments and discuss their results and describe our future study.
\section{Euler-Lagrange equation and a relaxation algorithm}

\smallskip

In this section we propose a method for solving the total variation minimization problem (\ref{original}) in 1-dimensional case. The details could be found in ~\cite{For07}. For gray-scale 1-dimensional inpainting problem, the functional (\ref{original}) becomes
\begin{equation}\label{energy}
\mathcal{J}(u):=\int_{\Omega}\left\vert 1_{\Omega\setminus D}(u(x)-g(x))\right\vert^{2}dx+2\lambda \int _{\Omega}\left\vert u'(x)\right\vert dx,
\end{equation}
where $D\subset\Omega$ is the damaged domain with measure $\mu(\Omega\setminus D)>0$, and $1_{\Omega\setminus D}$ denotes the characteristic function of $\Omega\setminus D$.

Associated to $\mathcal{J}$ we have the formal Euler-Lagrange equation:
\begin{equation}\label{Eulerform}
-\lambda (\frac{u'}{\left\vert u'\right\vert})'  + (u-g)1_{\Omega \backslash D}=0,
\end{equation}
with suitable boundary conditions. In our case, we use Neumann conditions.

Later we introduce a new functional given by
\begin{equation}\label{}
    \mathcal{\varepsilon}_{h}(u, w)=2\int _{\Omega}\left\vert 1_{\Omega\setminus D}(u(x)-g(x))\right\vert^{2}dx+2\lambda \int_{\Omega}(w\left\vert u'\right\vert^{2}+\frac{1}{w})dx,
\end{equation}
where $u\in W^{1,2}(\Omega; R)$, and $w\in L^{2}(\Omega; R)$ is such that $\epsilon_{h}\leq w \leq \frac{1}{\epsilon_{h}}$, where $\{\epsilon_{h}\}$ is a positive decreasing sequence such that $lim_{h\rightarrow \infty}\epsilon_{h}=0$. While the variable $u$ again is the function to be reconstructed, we call the variable $w$ the \verb"gradient weight".

For any given $u^{(0)}$ and $w^{(0)}$, we define the following iterative alternating-minimization algorithm:\\
\begin{equation}\label{algorithm}
 \left\lbrace
  \begin{array}{c l}
    u^{(n+1)}=arg\min _{u\in W^{1,2}(\Omega; R)}\varepsilon(u,w^{(n)}),\\
    \;\;\;w^{(n+1)}=arg \min_{\epsilon_{h}\leq w\leq \frac{1}{\epsilon_{h}}}\varepsilon(u^{(n+1)},w).
  \end{array}
\right. \end{equation}
\\
Then we have the 1-dimensional convergent result of Theorem 7.2 in ~\cite{For07}.
\begin{thm} \label{minimizer}
The sequence ${\{ {u^{(n)}}\} _{n \in \mathbb{N}}}$ has subsequences that converge strongly in ${L^2}(\Omega ;\;\mathbb{R})$
and weakly in ${W^{1,\;2}}(\Omega ;\;\mathbb{R})$ to a stationary point $u^{(\infty)}$ of $\mathcal{J}$; i.e., $u^{(\infty)}$ solves the Euler-Lagrange equations (\ref{Eulerform}). Moreover, if $\mathcal{J}$ has a unique minimizer $u^{*}$, then $u^{(\infty)}=u^{*}$ and the full sequence ${\{ {u^{(n)}}\} _{n \in \mathbb{N}}}$ converges to $u^{*}$.
\end{thm}

 From \textbf{Theorem \ref{minimizer}} we conclude that both $\mathcal{J}$ and $\mathcal{\varepsilon}_{h}(\cdot, w)$ admit minimizers, their uniqueness is equivalent to the uniqueness of the solutions of the corresponding Euler-Lagrange equation (\ref{Eulerform}). If uniqueness of the solution is satisfied, then the algorithm (\ref{algorithm})
can be reformulated equivalently as the following two-step iterative procedure:\\
\begin{itemize}
  \item Find $u^{(n+1)}$, which solves
\begin{equation}\label{modelproblem}
\int_{\Omega}( w^{(n)}(u^{(n+1)})'v'+\frac{1_{\Omega\backslash D}}{\lambda}(u^{(n+1)}-g)v)dx=0 \quad \forall v\in W^{1,2}(\Omega; R);
\end{equation}
  \item Compute directly $w^{(n+1)}$ by \\

\[ w^{(n+1)} =\epsilon_{h}\vee \frac{1}{\left\vert (u^{(n+1)})'\right\vert}\wedge \frac{1}{\epsilon_{h}}:= \left\lbrace
  \begin{array}{c l}
    \frac{1}{\left\vert (u^{(n+1)})'\right\vert}  & \text{if $\epsilon_{h} \leq\frac{1}{\left\vert (u^{(n+1)})'\right\vert}<\frac{1}{\epsilon_{h}}$},\\
    \epsilon_{h} & \text{if $\frac{1}{\left\vert (u^{(n+1)})'\right\vert}<\epsilon_{h}<\frac{1}{\epsilon_{h}}$},\\
    \frac{1}{\epsilon_{h}} & \text{otherwise}.
  \end{array}
\right. \]
\end{itemize}
In the following sections we illustrate the finite element approximation of the Euler-Lagrange equation (\ref{Eulerform}) similar to \cite[section 8]{For07}.
However, the interesting solutions may be discontinuous. In order to deal promptly with the discontinuity, we have studied the applications to gray-scale 1-dimensional inpainting problem by discontinuous Galerkin method for total variation minimization.

%Find $u^{(n+1)}$, which solves
%\begin{equation}\label{modelproblem}
%\int_{\Omega}( w^{(n)}(u^{(n+1)})'v'+\frac{1_{\Omega\backslash D}}{\lambda}(u^{(n+1)}-g)v)dx=0 \quad \forall v\in W^{1,2}(\Omega; R);
%\end{equation}
%Compute directly $w^{(n+1)}$ by \\
%
%\[ w^{(n+1)} =\epsilon_{h}\vee \frac{1}{\left\vert (u^{(n+1)})'\right\vert}\wedge \frac{1}{\epsilon_{h}}:= \left\lbrace
%  \begin{array}{c l}
%    \frac{1}{\left\vert (u^{(n+1)})'\right\vert}  & \text{if $\epsilon_{h} \leq\frac{1}{\left\vert (u^{(n+1)})'\right\vert}<\frac{1}{\epsilon_{h}}$},\\
%    \epsilon_{h} & \text{if $\frac{1}{\left\vert (u^{(n+1)})'\right\vert}<\epsilon_{h}<\frac{1}{\epsilon_{h}}$},\\
%    \frac{1}{\epsilon_{h}} & \text{otherwise}.
%  \end{array}
%\right. \]
%
%
%
%% -----------------------------------------------------------
%
\section{Finite element method for total variation minimization}
\subsection{Finite element method formulation for problem (\ref{Eulerform}).}
Denote $\widetilde{\lambda}=\frac{1_{\Omega\backslash D}}{\lambda}$, then for a given gradient weight $w^{(n)}$, the finite element method for solving  (\ref{modelproblem}) is to find $u^{(n+1)}$ such that
\begin{equation}\label{Femformulation}
    a(u^{(n+1)},v)=<F,v> \; \forall v\in W^{1,2}(\Omega; R),
\end{equation}
where $$a(u^{(n+1)}, v)= \int_{\Omega}( w^{(n)}(u^{(n+1)})'v'+\widetilde{\lambda}u^{(n+1)}v)dx$$ and $$<F, v>=\int_{\Omega}\widetilde{\lambda}gv dx.$$
Suppose the problem domain $\Omega$ is discretized into $N$ equal size of elements: \\ $0=x_{0}=x_{1}<\cdots<x_{N}=1$, denote $I_{m}=(x_{m},x_{m+1})$ and $h$ the mesh size. The integral for the $m^{th}$ element is
$$\int_{x_{m}}^{x_{m+1}}( w^{(n)}(u^{(n+1)})'v'+\widetilde{\lambda}u^{(n+1)}v)dx.$$
The trial function $u$ is expressed as
$$u^{(n+1)}=\phi_{1}u^{(n+1)}_{m}+\phi_{2}u^{(n+1)}_{m+1}$$
with the usual nodal basis functions
$$\phi_{1}(x)=\frac{x_{m+1}-x}{h}; \quad\phi_{2}(x)=\frac{x-x_{m}}{h}.$$
In our example (Section \ref{nume}), the value of $g$ in each element is a constant (we denote it by $\widetilde{g}$) and the value of $\widetilde{\lambda}$ is either $0$ or $\frac{1}{\lambda}$.
Now we could compute the element matrix and the element load vector
$$\mathbf{A}_{m}^{(n+1)} = w^{(n)}\left(
        \begin{array}{cc}
          \frac{1}{h} & -\frac{1}{h} \\
          -\frac{1}{h} & \frac{1}{h} \\
        \end{array}
      \right)+\widetilde{\lambda}\left(
                \begin{array}{cc}
                  \frac{h}{3} & \frac{h}{6} \\
                  \frac{h}{6} & \frac{h}{3} \\
                \end{array}
              \right); \quad \mathbf{b}_{m}^{(n+1)} = \widetilde{\lambda}\widetilde{g}\left(
                             \begin{array}{c}
                                \frac{h}{2} \\
                                \frac{h}{2} \\
                             \end{array}
                           \right).
$$
Assembling the element matrices and element load vectors, we could obtain the linear system
$\mathbf{A}^{(n+1)}\mathbf{u}^{(n+1)}=\mathbf{b}^{(n+1)}$.
\subsection{\textbf{Numerical implementation of the alternating-minimization algorithm.}}
\textbf{Input}: Data vector $\mathbf{\overline{g}}$, $\epsilon_{h}>0$, initial gradient weight $w^{(0)}$ with $\epsilon_{h}\leq w^{(0)}\leq \frac{1}{\epsilon_{h}}$, number $n_{max}$ of outer iterations.\\
\textbf{Parameters}: Positive weight $\tilde{\lambda}$.\\
\textbf{Output}: Approximation $u^{*}$ of the minimizer of ($\ref{energy}$).\\
$\mathbf{u}^{(0)}:=0$;\\
for $n:=0$ to $n_{max}$ do\\
Compute $\mathbf{u}^{(n+1)}$ such that $\mathbf{A}^{(n+1)}\mathbf{u}^{(n+1)}=\mathbf{b}$;\\
Compute the gradient ${(u^{(n+1)}|_{I_{m}})}'=u^{(n+1)}_{m}{\phi_{1}}'+u^{(n+1)}_{m+1}{\phi_{2}}'=-\frac{u^{(n+1)}_{m}}{h}+\frac{u^{(n+1)}_{m+1}}{h}$;\\
$w^{(n+1)}=\epsilon_{h}\vee \frac{1}{\left\vert (u^{(n+1)})'\right\vert}\wedge \frac{1}{\epsilon_{h}};$\\
endfor\\
$u^{*}:=u^{(n+1)}.$
\section{Discontinuous Galerkin method for total variation minimization}
\subsection{Discontinuous Galerkin method}
In this section, we will use Discontinuous Galerkin method to solve the same problem computing the solution of Euler-Lagrange equation (\ref{Eulerform}).
Let us consider the problem
\begin{equation}\label{classical}
    -(wu')'+\widetilde{\lambda}u=\widetilde{\lambda}g
\end{equation}
Let $0=x_{0}=x_{1}<\cdots<x_{N}=1$ be an uniform partition, denote $I_{n}=(x_{n},x_{n+1})$.
Denote by $\mathcal{D}_{1}$ the space of piecewise discontinuous polynomials of degree 1:
\begin{equation*}
\mathcal{D}_{1}=\{v:v|_{I_{n}}\in P_{1}(I_{n})\, \forall n=0, \ldots, N-1\},
\end{equation*}
where $ P_{1}(I_{n})$ is the space of polynomials of degree 1 on the interval $I_{n}$.\\

Then we define the jump and the average of $v$ at the boundary points of $I_{n}$:
\begin{equation*}
 [v(x_{n})]=v(x_{n}^{-})-v(x_{n}^{+}), \{v(x_{n})\}=\frac{1}{2}(v(x_{n}^{-})+v(x_{n}^{+}))\quad\forall n=1, \ldots, N-1.
\end{equation*}
We also extend the definition of jump and average at $x_{0}$ and $x_{N}$:
\begin{equation*}
 [v(x_{0})]=-v(x_{0}^{+}),\, \{v(x_{0})\}=v(x_{0}^{+}), \, [v(x_{N})]=v(x_{N}^{-}),\,\{v(x_{N})\}=v(x_{N}^{-}).
\end{equation*}

Next we introduce the penalty terms of the solution:
\begin{equation*}
J_{0}(u,v)=\sum_{n=0}^{N}\frac{\alpha}{h}[u(x_{n})][v(x_{n})].
\end{equation*}
where $\alpha$ is the real nonnegative number and $h$ is the mesh size.

Now we multiply (\ref{classical}) by $v\in\mathcal{D}_{1}$ and use integrating by parts on each interval $I_{n}$:
\begin{equation*}
    \int_{x_{n}}^{x_{n+1}}(wu'v'+\widetilde{\lambda}uv) dx-wu'v|_{x_{n}^{+}}^{x_{n+1}^{-}}= \int_{x_{n}}^{x_{n+1}}\widetilde{\lambda}gv.
\end{equation*}
By adding all $N$ equations above, we obtain
\begin{equation*}
     \sum_{n=0}^{N-1}\int_{x_{n}}^{x_{n+1}}(wu'v'+\widetilde{\lambda}uv) dx- \sum_{n=0}^{N-1}wu'v|_{x_{n}^{+}}^{x_{n+1}^{-}}= \int_{0}^{1}\widetilde{\lambda}gv.
\end{equation*}
Then we have
\begin{equation*}
     \sum_{n=0}^{N-1}\int_{x_{n}}^{x_{n+1}}(wu'v'+\widetilde{\lambda}uv) dx- \sum_{n=0}^{N}\{w(x_{n})u'(x_{n})\}[v(x_{n})]= \int_{0}^{1}\widetilde{\lambda}gv.
\end{equation*}
If $u$ is a solution of (\ref{classical}), then $u$ is continuous(\;$[u(x_{n})]=0$ for all $1\leq n \leq N-1$\;), thus $u$ satisfies
\begin{align*}
     &\sum_{n=0}^{N-1}\int_{x_{n}}^{x_{n+1}}(wu'v'+\widetilde{\lambda}uv) dx- \sum_{n=0}^{N}\{w(x_{n})u'(x_{n})\}[v(x_{n})]+\beta  \sum_{n=0}^{N}\{w(x_{n})v'(x_{n})\}[u(x_{n})]+J_{0}(u, v)\\&= \int_{0}^{1}\widetilde{\lambda}gv+\beta(-w(x_{0})v'(x_{0})u(x_{0})
     +w(x_{N})v'(x_{N})u(x_{N}))+\frac{\alpha}{h}(u(x_{0})v(x_{0})
     +u(x_{N})v(x_{N})).&
\end{align*}

Now the DG methods for solving (\ref{classical}) is to find $u\in \mathcal{D}_{1}$ such that
\begin{equation}\label{DGformulation}
    a(u,v)=<F,v> \; \forall v\in\mathcal{D}_{1},
\end{equation}
where
\begin{align*}
     a(u,v)=&\sum_{n=0}^{N-1}\int_{x_{n}}^{x_{n+1}}(wu'v'+\widetilde{\lambda}uv) dx- \sum_{n=0}^{N}\{w(x_{n})u'(x_{n})\}[v(x_{n})]\\&
     +\beta  \sum_{n=0}^{N}\{w(x_{n})v'(x_{n})\}[u(x_{n})]+J_{0}(u, v)&
\end{align*}
is the DG bilinear form, and
\begin{equation*}
    <F,\cdot>=\int_{0}^{1}\widetilde{\lambda}gv+\beta(-w(x_{0})v'(x_{0})u(x_{0})
     +w(x_{N})v'(x_{N})u(x_{N}))+\frac{\alpha}{h}(u(x_{0})v(x_{0})
     +u(x_{N})v(x_{N}))
\end{equation*}
is the linear form. In our example (Section \ref{nume}), we take the parameter $\beta=1$ so that the DG biliear form is symmetric.

\subsection{Linear system}
In this subsection, we derive the linear system obtained from the DG method. We choose for local basis functions of $P_{1}(I_{n})$ the nodal basis functions, i.e, $P_{1}(I_{n})=\text{span}\{\phi_{1}^{n},\phi_{2}^{n}\}$ with
\begin{equation*}
    \phi_{1}^{n}(x)=\frac{x_{i+1}-x}{x_{i+1}-x_{i}}; \quad \phi_{2}^{n}(x)=\frac{x-x_{i}}{x_{i+1}-x_{i}}.
\end{equation*}
The global basis functions $\{\Phi_{i}^{n}\}$ for the space $\mathcal{D}_{1}$ are obtained from the local basis functions by extending them by zero:
\[ \Phi_{i}^{n}(x) = \left\lbrace
  \begin{array}{c l}
     \phi_{i}^{n}(x) & \text{if $x\in I_{n}$},\\
    0 & \text{otherwise}.
  \end{array}
\right. \]
We can then expand the DG solution as
\begin{equation}\label{DGsolution}
    u^{DG}(x)=\sum_{m=0}^{N-1}\sum_{j=1}^{2}\alpha_{j}^{m}\Phi_{j}^{m}(x).
\end{equation}
Inserting this form of $u^{DG}$ into the scheme (\ref{DGformulation}), we get
\begin{equation*}
    \sum_{m=0}^{N-1}\sum_{j=1}^{2}\alpha_{j}^{m}a(\Phi_{j}^{m},\Phi_{i}^{n})=<F,\Phi_{i}^{n}> , \; \forall \;0\leq n \leq N-1,  \forall\; 1\leq i \leq 2.
\end{equation*}
We would then obtain a linear system $\mathbf{A\alpha}=\mathbf{b}$, where $\mathbf{\alpha}$ is the vector with components $\alpha_{j}^{m}$, $\mathbf{A}$ is the matrix with entries $a(\Phi_{j}^{m},\Phi_{i}^{n})$, and $\mathbf{b}$ is the vector with the components $<F,\Phi_{i}^{n}>$.

\subsubsection{\textbf{Computing the matrix} $\mathbf{A}$}
In this section, we will first show how to compute the local matrices. We will regroup the terms $a(\Phi_{j}^{m},\Phi_{i}^{n})$ into three groups: the terms involving integrals over $I_{n}$, the terms involving the interior nodes $x_{n}$, and the terms involving the boundary nodes $x_{0}$ and $x_{N}$.

Firstly, we consider the term corresponding to the integrals over $I_{n}$. On each element $I_{n}$, the DG solution $u^{DG}$ can be expressed as
\begin{equation}\label{localsolution}
    u^{DG}(x)=\alpha_{1}^{n}\phi_{1}^{n}(x)+\alpha_{2}^{n}\phi_{2}^{n}(x)\;\forall x\in I_{n}.
\end{equation}
Thus, using (\ref{localsolution}) and choosing $v=\phi_{i}^{n}$ for $i=1,2$, we get
\begin{equation*}
\int_{x_{n}}^{x_{n+1}}(w{(u^{DG})}'{(\phi_{i}^{n})}'+\widetilde{\lambda}u^{DG}\phi_{i}^{n}) dx =\sum_{j=1}^{2}\alpha_{j}^{n}\int_{x_{n}}^{x_{n+1}}(w{(\phi_{j}^{n})}'{(\phi_{i}^{n})}'+\widetilde{\lambda}\phi_{j}^{n}\phi_{i}^{n}) dx \;\forall i=1,2.
\end{equation*}
This linear system can be written as $\mathbf{A}_{n}\mathbf{\alpha}^{n}$, where
\begin{equation*}
    (\mathbf{A}_{n})_{ij}=\int_{x_{n}}^{x_{n+1}}(w{(\phi_{j}^{n})}'{(\phi_{i}^{n})}'+\widetilde{\lambda}\phi_{j}^{n}\phi_{i}^{n}) dx,\quad \mathbf{\alpha}^{n}=\left(
                                   \begin{array}{c}
                                     \alpha_{1}^{n} \\
                                     \alpha_{2}^{n} \\
                                   \end{array}
                                 \right).
\end{equation*}
We could compute the $\mathbf{A}_{n}$:
\begin{equation*}
    \mathbf{A}_{n}=\frac{w_{n}}{h}\left(
                                    \begin{array}{cc}
                                      1 & -1 \\
                                      -1 & 1 \\
                                    \end{array}
                                  \right)+\tilde{\lambda}h\left(
                                                           \begin{array}{cc}
                                                             1/3 & 1/6 \\
                                                             1/6 & 1/3 \\
                                                           \end{array}
                                                         \right).
\end{equation*}

Second, we consider the terms involving the interior nodes $x_{n}$. Let us express
\begin{align*}
    &-\{w(x_{n})(u^{DG})'(x_{n})\}[v(x_{n})]+\beta\{w(x_{n})v'(x_{n})\}[u^{DG}(x_{n})]+\frac{\alpha}{h}[u^{DG}(x_{n})][v(x_{n})]\\
    &=b_{n}+c_{n}+d_{n}+e_{n},
\end{align*}
where the terms are defined as below:
\begin{eqnarray*}
  b_{n} = \frac{1}{2}w(x_{n}^{+})(u^{DG})'(x_{n}^{+})v(x_{n}^{+})-\frac{\beta}{2}w(x_{n}^{+})u^{DG}(x_{n}^{+})v'(x_{n}^{+})+\frac{\alpha}{h}u^{DG}(x_{n}^{+})v'(x_{n}^{+}), \\
  c_{n} = -  \frac{1}{2}w(x_{n}^{-})(u^{DG})'(x_{n}^{-})v(x_{n}^{-})+\frac{\beta}{2}w(x_{n}^{-})u^{DG}(x_{n}^{-})v'(x_{n}^{-})+\frac{\alpha}{h}u^{DG}(x_{n}^{-})v'(x_{n}^{-}) , \\
  d_{n} = -  \frac{1}{2}w(x_{n}^{+})(u^{DG})'(x_{n}^{+})v(x_{n}^{-})-\frac{\beta}{2}w(x_{n}^{-})u^{DG}(x_{n}^{+})v'(x_{n}^{-})-\frac{\alpha}{h}u^{DG}(x_{n}^{+})v'(x_{n}^{-}), \\
  e_{n} =   \frac{1}{2}w(x_{n}^{-})(u^{DG})'(x_{n}^{-})v(x_{n}^{+})+\frac{\beta}{2}w(x_{n}^{+})u^{DG}(x_{n}^{-})v'(x_{n}^{+})-\frac{\alpha}{h}u^{DG}(x_{n}^{-})v'(x_{n}^{+}) .
\end{eqnarray*}
Now with the expression (\ref{localsolution}) and the choice $v=\phi_{i}^{n}$, the four terms defined above yields the local matrices $\mathbf{B}_{n}, \mathbf{C}_{n}, \mathbf{D}_{n}$ and $\mathbf{E}_{n}$ which are:
\begin{eqnarray*}
  \mathbf{B}_{n} = \frac{1}{2h}\left(
                                 \begin{array}{cc}
                                   -w(x_{n}^{+})+\beta w(x_{n}^{+})+2\alpha & w(x_{n}^{+}) \\
                                   -\beta w(x_{n}^{+}) & 0 \\
                                 \end{array}
                               \right),
   \\
     \mathbf{C}_{n} = \frac{1}{2h}\left(
                                 \begin{array}{cc}
                                   0 & -\beta w(x_{n}^{-}) \\
                                    w(x_{n}^{-}) & -w(x_{n}^{+})+\beta w(x_{n}^{+})+2\alpha  \\
                                 \end{array}
                               \right), \\
  \mathbf{D}_{n} = \frac{1}{2h}\left(
                                 \begin{array}{cc}
                                    \beta w(x_{n}^{-}) & 0 \\
                                  w(x_{n}^{+})-\beta w(x_{n}^{-})-2\alpha & -w(x_{n}^{+}) \\
                                 \end{array}
                               \right) ,\\
  \mathbf{E}_{n} = \frac{1}{2h}\left(
                                 \begin{array}{cc}
                                   -w(x_{n}^{-}) &  w(x_{n}^{-})-\beta w(x_{n}^{+})-2\alpha\\
                                    0 & \beta w(x_{n}^{+}) \\
                                 \end{array}
                               \right).
\end{eqnarray*}

Finally, we compute the local matrices from the boundary nodes $x_{0}$ and $x_{N}$:
\begin{eqnarray*}
  f_{0} &=&   w(x_{0})(u^{DG})'(x_{0})v(x_{0})-\beta w(x_{0})u^{DG}(x_{0})v'(x_{0})+\frac{\alpha}{h}u^{DG}(x_{0})v'(x_{0}),  \\
  f_{N} &=&   -w(x_{N})(u^{DG})'(x_{N})v(x_{N})+\beta w(x_{N})(u^{DG})(x_{N})v'(x_{N})+\frac{\alpha}{h}u^{DG}(x_{N})v'(x_{N}),
\end{eqnarray*}
which yields the local matrices $\mathbf{F}_{0}$ and $\mathbf{F}_{N}$:
\begin{eqnarray*}
  \mathbf{F}_{0} = \frac{1}{h}\left(
                                 \begin{array}{cc}
                                   -w(x_{0})+\beta w(x_{0})+\alpha & w(x_{0}) \\
                                   -\beta w(x_{0}) & 0 \\
                                 \end{array}
                               \right),
   \\
     \mathbf{F}_{N} = \frac{1}{h}\left(
                                 \begin{array}{cc}
                                   0 & -\beta w(x_{N}) \\
                                    w(x_{N}) & -w(x_{N})+\beta w(x_{N})+\alpha  \\
                                 \end{array}
                               \right). \\
\end{eqnarray*}

Assuming that the unknowns are listed in the following order:
\begin{equation*}
    (\alpha_{1}^{0},\alpha_{2}^{0},\alpha_{1}^{1},\alpha_{2}^{1},\alpha_{1}^{2},\alpha_{2}^{2},\ldots,\alpha_{1}^{N-1},\alpha_{2}^{N-1}),
\end{equation*}
we obtain the global matrix $\mathbf{A}$ which is block tridiagonal
\begin{equation*}
\mathbf{A}=\left(
  \begin{array}{cccccc}
    \mathbf{M}_{0} & \mathbf{D}_{1} & {} & {} & {} & {} \\
    \mathbf{E}_{1} & \mathbf{M} & \mathbf{D}_{2} & {}& {} &{}  \\
     & \cdots & \cdots & \cdots & {} &{}  \\
     & {} & \cdots & \cdots & \cdots &{}  \\
     & {} &{}  & \mathbf{E}_{N-2} & \mathbf{M} & \mathbf{D}_{N-1} \\
     & {} &{}  & {} & \mathbf{E}_{N-1} & \mathbf{M}_{N} \\
  \end{array}
\right),
\end{equation*}
where
\begin{equation*}
    \mathbf{M}=\mathbf{A}_{n}+\mathbf{B}_{n}+\mathbf{C}_{n+1}, \mathbf{M}_{0}=\mathbf{A}_{0}+\mathbf{F}_{0}+\mathbf{C}_{1},
    \mathbf{M}_{N}=\mathbf{A}_{N-1}+\mathbf{F}_{N}+\mathbf{C}_{N-1}.
\end{equation*}

\subsubsection{\textbf{Computing the right hand side} $\mathbf{b}$}
% -----------------------------------------------------------
Each component of $\mathbf{b}$ can be obtained by computing
\begin{align*}
    &<F,\Phi_{n}^{i}>=\int_{0}^{1}\widetilde{\lambda}g\Phi_{n}^{i}dx+\beta(-w(x_{0})(\Phi_{n}^{i})'(x_{0})u(x_{0})
     +w(x_{N})(\Phi_{n}^{i})'(x_{N})u(x_{N}))\\&\qquad\qquad\qquad+\frac{\alpha}{h}(u(x_{0})\Phi_{n}^{i}(x_{0})
     +u(x_{N})\Phi_{n}^{i}(x_{N})).
\end{align*}
Because of the local support of $\Phi_{n}^{i}$, the first term is reduced to
\begin{equation*}
    \int_{0}^{1}\widetilde{\lambda}g\Phi_{n}^{i}dx=\int_{x_{n}}^{x_{n+1}}\widetilde{\lambda}g\phi_{n}^{i}dx.
\end{equation*}
We arrange the components of $\mathbf{b}$ in an order consistent with the order of the unknowns $\alpha_{i}^{n}$:
\begin{equation*}
    (b_{1}^{0},b_{2}^{0},b_{1}^{1},b_{2}^{1},b_{1}^{2},b_{2}^{2},\ldots,b_{1}^{N-1},b_{2}^{N-1}),
\end{equation*}
where the first two components and last two components are
\begin{eqnarray*}
  b_{1}^{0}=\frac{\tilde{\lambda}g_{0}h}{2}+\frac{\beta w(x_{0})u(x_{0})}{h}+\frac{\alpha u(x_{0})}{h}, \\
  b_{2}^{0}=\frac{\tilde{\lambda}g_{0}h}{2}-\frac{\beta w(x_{0})u(x_{0})}{h},\\
  b_{1}^{N-1}=\frac{\tilde{\lambda}g_{N-1}h}{2}-\frac{\beta w(x_{N})u(x_{N})}{h}, \\
  b_{2}^{N-1}=\frac{\tilde{\lambda}g_{N-1}h}{2}+\frac{\beta w(x_{N})u(x_{N})}{h}+\frac{\alpha u(x_{N})}{h},
\end{eqnarray*}
and the other $2(N-2)$ components are
\begin{equation*}
    b_{i}^{n}=\frac{\tilde{\lambda}g_{n}h}{2}, \quad \forall 1\leq n\leq N-2, \forall 1\leq i\leq 2.
\end{equation*}
\section{Numerical examples and future study}$\label{nume}$
In this section we show numerical results of the applications of the finite element method and Discontinuous Galerkin method for total variation minimization described above.

We would first use finite element method for total variation minimization the do some experiments with different values of parameter $\widetilde{\lambda}$ in the alternating minimization algorithm. Let us consider the same signal with the same inpainting interval $(\frac{1}{3},\frac{2}{3})$, see Figure $\ref{nume}.1$. We conclude that quality of the recovery image is becoming better as the value of $\widetilde{\lambda}$ increases.
\begin{center}
  \begin{tabular}{ccc}
    \includegraphics[width=1.6in]{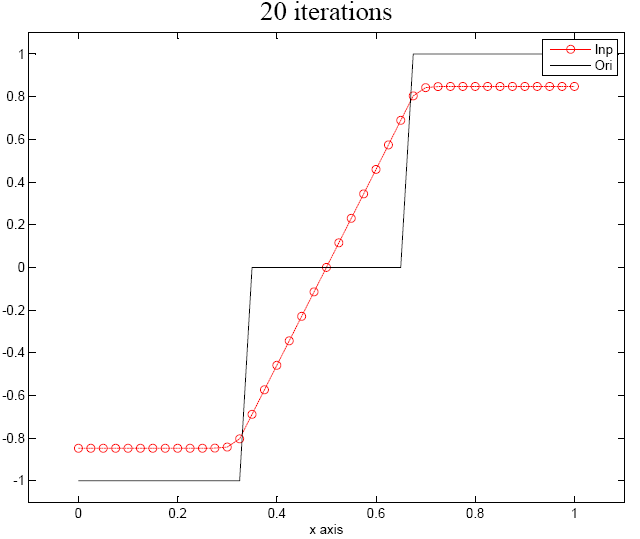} &
    \includegraphics[width=1.6in]{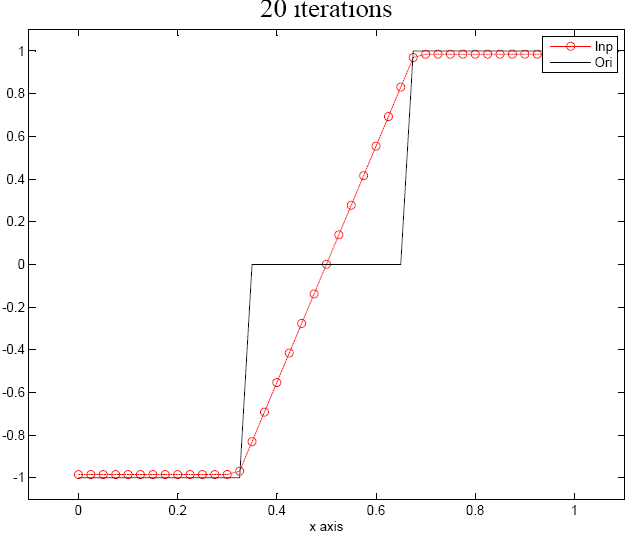} &
    \includegraphics[width=1.6in]{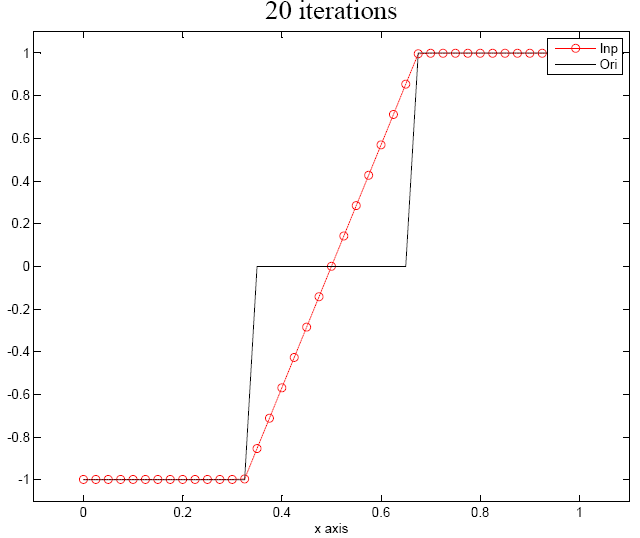}
  \end{tabular}
\end{center}
\bigskip

\noindent Figure $\ref{nume}.1$ \emph{The three pictures from left to right represent the result of three different values of $\widetilde{\lambda}$ ($\widetilde{\lambda}=10, 100, 1000$), respectively.}

\bigskip

Second, we fix the outer iteration to 20 and compare the convergence speed for three different values of $\widetilde{\lambda}$ ($\widetilde{\lambda}=10, 100, 1000$), see the left picture of Figure $\ref{nume}.2$. And we conclude that the convergence speed increases as the value of $\widetilde{\lambda}$ increases.

Finally, we modify $w^{(n+1)}$ from $\epsilon_{h}\vee \frac{1}{\left\vert (u^{(n+1)})'\right\vert}\wedge \frac{1}{\epsilon_{h}}$ to $(\epsilon_{h}\vee \frac{1}{\left\vert (u^{(n+1)})'\right\vert}\wedge
\frac{1}{\epsilon_{h}})^{2-\tau}$, and check the convergence speed for different value of $\tau$($\tau$=1,\;0.9,\;0.8,\;0.7,\;0.6,\;0.5). The result is shown in the right picture of Figure $\ref{nume}.2$. And we conclude that the convergence speed increases as the value of $\tau$ decreases.

\begin{center}
  \begin{tabular}{cc}
    \includegraphics[width=2.2in]{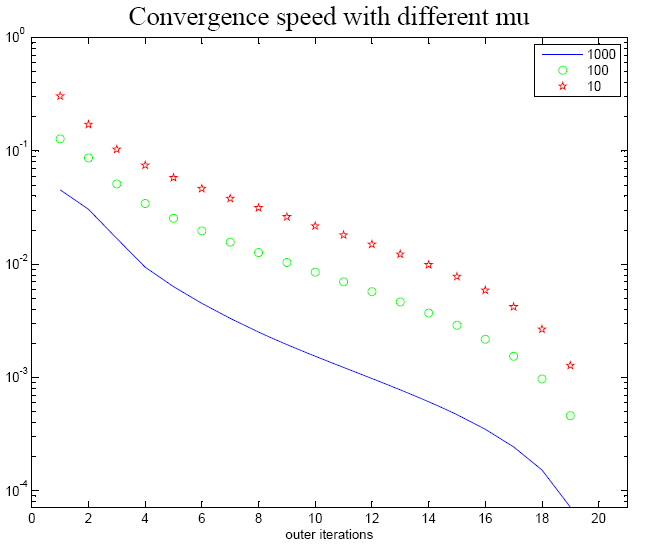} &
    \includegraphics[width=2.2in]{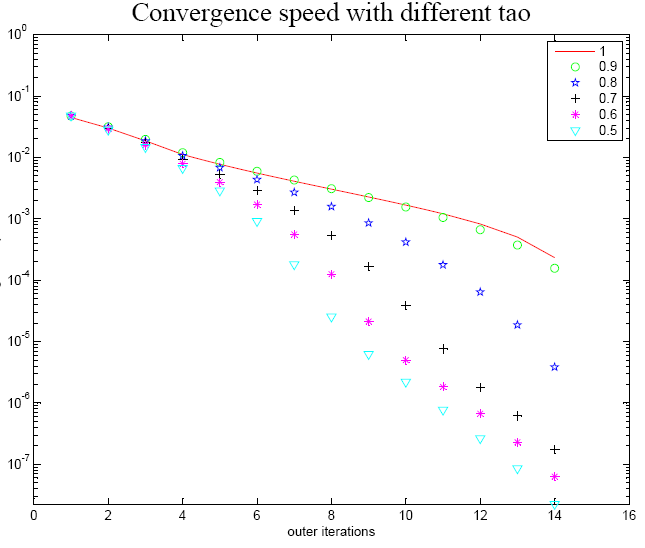}
  \end{tabular}
\end{center}

\bigskip

\noindent Figure $\ref{nume}.2$ \emph{The left picture illustrates the convergence speed for three different values of $\widetilde{\lambda}$ ($\widetilde{\lambda}=10, 100, 1000$); the right picture shows the convergence speed for different value of $\tau$ ($\tau$=1,\;0.9,\;0.8,\;0.7,\;0.6,\;0.5).}
\bigskip

Next we consider a signal with a jump (see Figure $\ref{nume}.3$).

\begin{center}
  \includegraphics[width=2.25in]{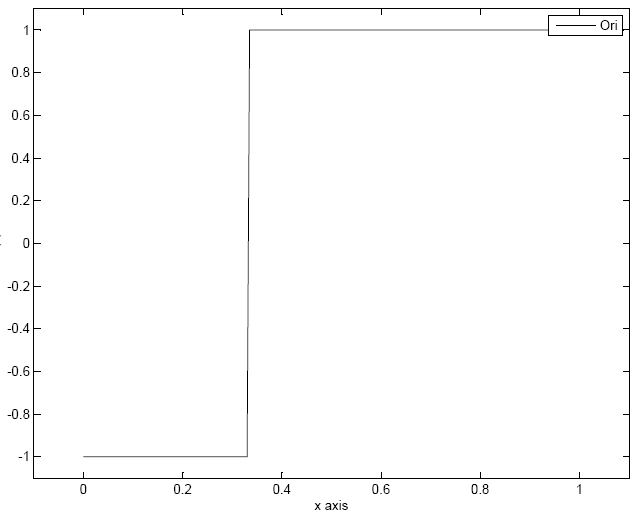}\\
\end{center}
\bigskip
\noindent Figure $\ref{nume}.3$ \emph{The signal of a step function.}
\bigskip

Let us recover the signal(Figure $\ref{nume}.3$) by finite element method and discontinuous Galerkin method for total variation, respectively. The results are shown in  Figure $\ref{nume}.4$. We observe that the finite element method for total variation minimization couldn't preserve the jump very well from our example. However, the discontinuous Galerkin method for total variation minimization preserves the jump rather well.

Our future study aims at the construction, analysis and implementation of new adaptive discontinuous Galerkin (DG) solvers for total variation minimization problems in two space dimensions. These methods are based on re-weighted least squares and are implemented by nester outer and inner iterations. The adaptivity concerns not only the space discretization but also the parameters involved in the inner and outer iterations as well as in the DG discretization. The inner iteration should be robust with respect to both the discretization and the gradient weights, and the number of inner iterations should be controlled in proper way. The robustness property is obviously connected with the preconditioning.
\begin{center}
  \begin{tabular}{cc}
    \includegraphics[width=2.2in]{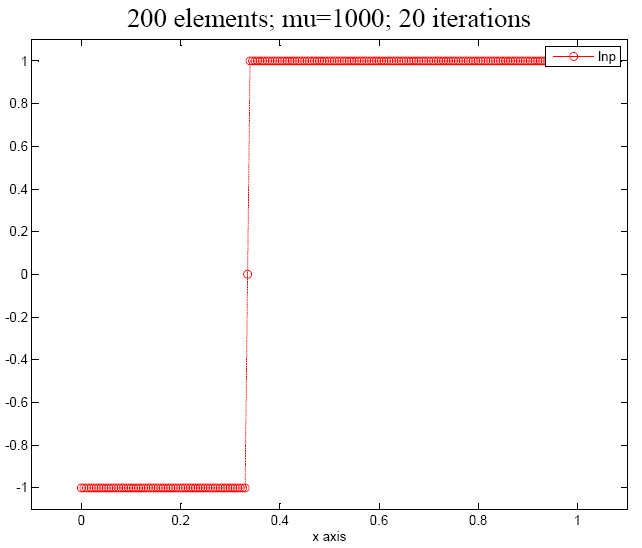} &
    \includegraphics[width=2.2in]{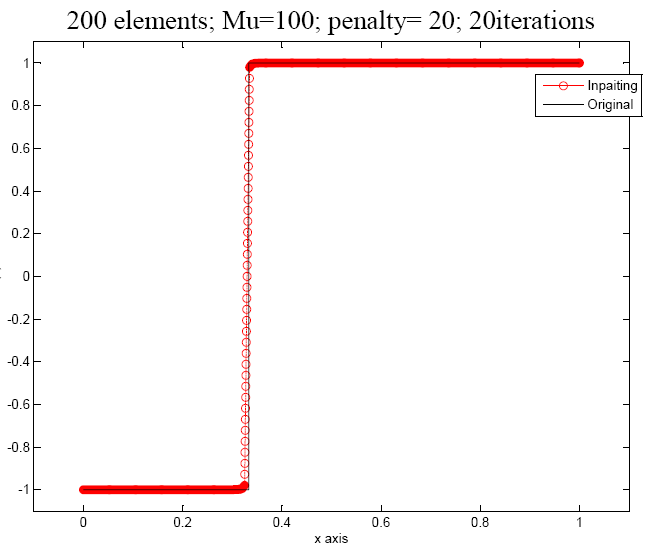}
  \end{tabular}
\end{center}

\bigskip

\noindent Figure $\ref{nume}.4$ \emph{The left picture shows the recovery image by the finite element method for total variation minimization; and the right picture illustrates discontinuous Galerkin method for total variation minimization. }
\bigskip

\section*{Acknowledgment}
The author wishes to thank Dr. Massimo Fornasier and Prof. Dr. Ulrich Langer for their valuable suggestions during the Project Seminar on Numerical Analysis in Johannes Kepler Linz University. (See also http://www.numa.uni-linz.ac.at/Teaching/LVA/2009s/Projsem/)

%% References

\bibliographystyle{model1-num-names}
\bibliography{Xbib}






% \end{thebibliography}

\end{document}